\documentclass{gtart}

\def\ifplaintex{\expandafter\ifx\csname documentclass\endcsname\relax}

\def\gtp{{\mathsurround=0pt\it $\cal G\mskip-2mu$eometry \&\ 
$\cal T\!\!$opology $\cal P\!$ublications}}  

\def\Addressesr{\bigskip
{\small \parskip 0pt \leftskip 0pt \rightskip 0pt plus 1fil \def\\{\par}
\sl\theaddress\par
\medskip
\rm Email:\stdspace\tt\theemail\hfill\rm Received:\qua\receiveddate \par}}

\def\recd{{\small Received:\qua\receiveddate\ifx\reviseddate\relax
\else\qquad Revised:\qua\reviseddate\fi\par}} 

\def\lognumber#1{\def\thelognumber{#1}}
\def\volumenumber#1{\def\thevolumenumber{#1}}
\def\volumeyear#1{\def\thevolumeyear{#1}}
\def\papernumber#1{\def\thepapernumber{#1}}
\def\pagenumbers#1#2{\def\startpage{#1}\def\finishpage{#2}}
\def\published#1{\def\publishdate{#1}}

\def\received#1{\def\receiveddate{#1}}

\def\accepted#1{\def\accepteddate{#1}}

\long\def\asciiabstract#1{\long\def\theasciiabstract{#1}}

\let\\\par\let\thelognumber\relax\let\thevolumenumber\relax
\let\thepapernumber\relax\let\thevolumeyear\relax\let\startpage\relax
\let\finishpage\relax\let\publishdate\relax\let\receiveddate\relax
\let\reviseddate\relax\let\accepteddate\relax\let\theasciititle\relax
\let\theasciiauthors\relax
\let\theasciiabstract\relax

\let\theasciiemail\relax

\ifplaintex
\font\logobig=cmssbx10 scaled 3836
\font\logomed=cmssbx10 scaled 2557
\else
\font\logobig=cmssbx10 scaled 4200
\font\logomed=cmssbx10 scaled 2800
\fi

\long\def\makeagttitle{   
\count0=\startpage
\agt\hfill      
\hbox to 45truept{\vbox to 0pt{\vglue -13truept{\logomed A\kern -.37em{\logobig 
T}\kern -.38em G}\vss}\hss}
\break
{\small Volume \thevolumenumber\ (\thevolumeyear)
\startpage--\finishpage\nl
Published: \publishdate}

\vglue .25truein

{\parskip=0pt\leftskip 0pt plus
1fil\def\\{\par\smallskip}{\Large\bf\thetitle}\par\medskip} \vglue
0.05truein

{\parskip=0pt\leftskip 0pt plus 1fil\def\\{\par}{\sc\theauthors}
\par\medskip}%
 
\vglue 0.03truein 

{\small\leftskip 25truept\rightskip 25truept{\bf Abstract}\stdspace\theabstract

{\bf AMS Classification}\stdspace\theprimaryclass
\ifx\thesecondaryclass\relax\else; \thesecondaryclass\fi\par
{\bf Keywords}\stdspace \thekeywords\par}\vglue 7truept

}   

\ifplaintex
\font\phead=cmsl9 scaled 950
\font\pnum=cmbx10 scaled 913
\font\pfoot=cmsl9 scaled 950
\headline{\vbox to 0pt{\vskip -4.5mm\line{\small\phead\ifnum
\count0=\startpage ISSN 1472-2739 (on-line) 1472-2747 (printed)
\hfill {\pnum\folio}\else\ifodd\count0\def\\{ }%
\ifx\theshorttitle\relax\thetitle\else\theshorttitle\fi\hfill{\pnum\folio}
\else\def\\{ and }{\pnum\folio}\hfill\ifx\theshortauthors\relax\theauthors
\else\theshortauthors\fi\fi\fi}\vss}}
\footline{\vbox to 0pt{\vglue 0mm\line{\small\pfoot\ifnum\count0=\startpage
\copyright\ \gtp\hfill\else
\agt, Volume \thevolumenumber\ (\thevolumeyear)\hfill\fi}\vss}}
\else
\font=cmsl9 scaled 1050
\font\lnum=cmbx10 
\font=cmsl9 scaled 1050
\makeatletter
\def\@oddhead{{\small\lhead\ifnum\count0=\startpage ISSN 1472-2739 
(on-line) 1472-2747 (printed)\hfill {\lnum\number\count0}\else\ifodd\count0
\def\\{ }\ifx\theshorttitle\relax \thetitle \else\theshorttitle\fi\hfill
{\lnum\number\count0}\else\def\\{ and }{\lnum\number\count0}
\hfill\ifx\theshortauthors\relax 
\theauthors\else\theshortauthors\fi\fi\fi}}\def\@evenhead{\@oddhead}
\def\@oddfoot{\small\lfoot\ifnum\count0=\startpage\copyright\ \gtp\hfill\else
\agt, Volume \thevolumenumber\ (\thevolumeyear)\hfill\fi}
\def\@evenfoot{\@oddfoot}
\makeatother
\fi
\let\maketitlepage\makeagttitle
\let\makeshorttitle\maketitlepage
\let\maketitle\maketitlepage

\newwrite\gtoutfile
\long\gdef\makeheadfile{  
{\def\\{, }\def\s{ }
\immediate\openout\gtoutfile head.xxx
\immediate\write\gtoutfile{To: math@arxiv.org}
\immediate\write\gtoutfile{Subject: put OR rep NNNNN:ppppp}
\immediate\write\gtoutfile{--text follows this line--}
\immediate\write\gtoutfile{Proxy-for: \ifx\theasciiauthors\relax
\theauthors\else\theasciiauthors\fi\s<\ifx\theasciiemail\relax\theemail\else\theasciiemail\fi>}
\immediate\write\gtoutfile{\noexpand\\}
\immediate\write\gtoutfile{Authors: \ifx\theasciiauthors\relax
\theauthors\else\theasciiauthors\fi}
{\def\\{ }\immediate\write\gtoutfile{Title: \ifx\theasciititle\relax
\thetitle\else\theasciititle\fi}}
\immediate\write\gtoutfile{Subj-class: GT or SG, GR etc}
\immediate\write\gtoutfile{MSC-class: \theprimaryclass\ifx\thesecondaryclass\relax\else, \thesecondaryclass\fi}
\immediate\write\gtoutfile{Journal-ref: Algebr. Geom. Topol. \thevolumenumber\s
(\thevolumeyear) \startpage-\finishpage}
\immediate\write\gtoutfile{Comments: Published by Algebraic and
Geometric Topology at}
\immediate\write\gtoutfile{\s\s\s  http://www.maths.warwick.ac.uk/agt/AGTVol\thevolumenumber/agt-\thevolumenumber-\thepapernumber.abs.html}
\immediate\write\gtoutfile{\noexpand\\}
\immediate\write\gtoutfile{}
\ifx\theasciiabstract\relax
\immediate\write\gtoutfile{\theabstract}\else
\immediate\write\gtoutfile{\theasciiabstract}\fi
\immediate\write\gtoutfile{}
\immediate\write\gtoutfile{\noexpand\\}
\immediate\write\gtoutfile{}
\immediate\closeout\gtoutfile}}  

\def\maketitlepage{\makeagttitle\makeheadfile}
\let\makeshorttitle\maketitlepage
\let\maketitle\maketitlepage

\def\ifplaintex{\expandafter\ifx\csname documentclass\endcsname\relax}

\def\gtp{{\mathsurround=0pt\it $\cal G\mskip-2mu$eometry \&\ 
$\cal T\!\!$opology $\cal P\!$ublications}}  

\def\Addressesr{\bigskip
{\small \parskip 0pt \leftskip 0pt \rightskip 0pt plus 1fil \def\\{\par}
\sl\theaddress\par
\medskip
\rm Email:\stdspace\tt\theemail\hfill\rm Received:\qua\receiveddate \par}}

\def\recd{{\small Received:\qua\receiveddate\ifx\reviseddate\relax
\else\qquad Revised:\qua\reviseddate\fi\par}} 

\def\lognumber#1{\def\thelognumber{#1}}
\def\volumenumber#1{\def\thevolumenumber{#1}}
\def\volumeyear#1{\def\thevolumeyear{#1}}
\def\papernumber#1{\def\thepapernumber{#1}}
\def\pagenumbers#1#2{\def\startpage{#1}\def\finishpage{#2}}
\def\published#1{\def\publishdate{#1}}

\def\received#1{\def\receiveddate{#1}}

\def\accepted#1{\def\accepteddate{#1}}

\long\def\asciiabstract#1{\long\def\theasciiabstract{#1}}

\let\\\par\let\thelognumber\relax\let\thevolumenumber\relax
\let\thepapernumber\relax\let\thevolumeyear\relax\let\startpage\relax
\let\finishpage\relax\let\publishdate\relax\let\receiveddate\relax
\let\reviseddate\relax\let\accepteddate\relax\let\theasciititle\relax
\let\theasciiauthors\relax
\let\theasciiabstract\relax

\let\theasciiemail\relax

\ifplaintex
\font\logobig=cmssbx10 scaled 3836
\font\logomed=cmssbx10 scaled 2557
\else
\font\logobig=cmssbx10 scaled 4200
\font\logomed=cmssbx10 scaled 2800
\fi

\long\def\makeagttitle{   
\count0=\startpage
\agt\hfill      
\hbox to 45truept{\vbox to 0pt{\vglue -13truept{\logomed A\kern -.37em{\logobig 
T}\kern -.38em G}\vss}\hss}
\break
{\small Volume \thevolumenumber\ (\thevolumeyear)
\startpage--\finishpage\nl
Published: \publishdate}

\vglue .25truein

{\parskip=0pt\leftskip 0pt plus
1fil\def\\{\par\smallskip}{\Large\bf\thetitle}\par\medskip} \vglue
0.05truein

{\parskip=0pt\leftskip 0pt plus 1fil\def\\{\par}{\sc\theauthors}
\par\medskip}%
 
\vglue 0.03truein 

{\small\leftskip 25truept\rightskip 25truept{\bf Abstract}\stdspace\theabstract

{\bf AMS Classification}\stdspace\theprimaryclass
\ifx\thesecondaryclass\relax\else; \thesecondaryclass\fi\par
{\bf Keywords}\stdspace \thekeywords\par}\vglue 7truept

}   

\ifplaintex
\font\phead=cmsl9 scaled 950
\font\pnum=cmbx10 scaled 913
\font\pfoot=cmsl9 scaled 950
\headline{\vbox to 0pt{\vskip -4.5mm\line{\small\phead\ifnum
\count0=\startpage ISSN 1472-2739 (on-line) 1472-2747 (printed)
\hfill {\pnum\folio}\else\ifodd\count0\def\\{ }%
\ifx\theshorttitle\relax\thetitle\else\theshorttitle\fi\hfill{\pnum\folio}
\else\def\\{ and }{\pnum\folio}\hfill\ifx\theshortauthors\relax\theauthors
\else\theshortauthors\fi\fi\fi}\vss}}
\footline{\vbox to 0pt{\vglue 0mm\line{\small\pfoot\ifnum\count0=\startpage
\copyright\ \gtp\hfill\else
\agt, Volume \thevolumenumber\ (\thevolumeyear)\hfill\fi}\vss}}
\else
\font=cmsl9 scaled 1050
\font\lnum=cmbx10 
\font=cmsl9 scaled 1050
\makeatletter
\def\@oddhead{{\small\lhead\ifnum\count0=\startpage ISSN 1472-2739 
(on-line) 1472-2747 (printed)\hfill {\lnum\number\count0}\else\ifodd\count0
\def\\{ }\ifx\theshorttitle\relax \thetitle \else\theshorttitle\fi\hfill
{\lnum\number\count0}\else\def\\{ and }{\lnum\number\count0}
\hfill\ifx\theshortauthors\relax 
\theauthors\else\theshortauthors\fi\fi\fi}}\def\@evenhead{\@oddhead}
\def\@oddfoot{\small\lfoot\ifnum\count0=\startpage\copyright\ \gtp\hfill\else
\agt, Volume \thevolumenumber\ (\thevolumeyear)\hfill\fi}
\def\@evenfoot{\@oddfoot}
\makeatother
\fi
\let\maketitlepage\makeagttitle
\let\makeshorttitle\maketitlepage
\let\maketitle\maketitlepage

\newwrite\gtoutfile
\long\gdef\makeheadfile{  
{\def\\{, }\def\s{ }
\immediate\openout\gtoutfile head.xxx
\immediate\write\gtoutfile{To: math@arxiv.org}
\immediate\write\gtoutfile{Subject: put OR rep NNNNN:ppppp}
\immediate\write\gtoutfile{--text follows this line--}
\immediate\write\gtoutfile{Proxy-for: \ifx\theasciiauthors\relax
\theauthors\else\theasciiauthors\fi\s<\ifx\theasciiemail\relax\theemail\else\theasciiemail\fi>}
\immediate\write\gtoutfile{\noexpand\\}
\immediate\write\gtoutfile{Authors: \ifx\theasciiauthors\relax
\theauthors\else\theasciiauthors\fi}
{\def\\{ }\immediate\write\gtoutfile{Title: \ifx\theasciititle\relax
\thetitle\else\theasciititle\fi}}
\immediate\write\gtoutfile{Subj-class: GT or SG, GR etc}
\immediate\write\gtoutfile{MSC-class: \theprimaryclass\ifx\thesecondaryclass\relax\else, \thesecondaryclass\fi}
\immediate\write\gtoutfile{Journal-ref: Algebr. Geom. Topol. \thevolumenumber\s
(\thevolumeyear) \startpage-\finishpage}
\immediate\write\gtoutfile{Comments: Published by Algebraic and
Geometric Topology at}
\immediate\write\gtoutfile{\s\s\s  http://www.maths.warwick.ac.uk/agt/AGTVol\thevolumenumber/agt-\thevolumenumber-\thepapernumber.abs.html}
\immediate\write\gtoutfile{\noexpand\\}
\immediate\write\gtoutfile{}
\ifx\theasciiabstract\relax
\immediate\write\gtoutfile{\theabstract}\else
\immediate\write\gtoutfile{\theasciiabstract}\fi
\immediate\write\gtoutfile{}
\immediate\write\gtoutfile{\noexpand\\}
\immediate\write\gtoutfile{}
\immediate\closeout\gtoutfile}}  

\def\maketitlepage{\makeagttitle\makeheadfile}
\let\makeshorttitle\maketitlepage
\let\maketitle\maketitlepage

\def\ifplaintex{\expandafter\ifx\csname documentclass\endcsname\relax}

\def\gtp{{\mathsurround=0pt\it $\cal G\mskip-2mu$eometry \&\ 
$\cal T\!\!$opology $\cal P\!$ublications}}  

\def\Addressesr{\bigskip
{\small \parskip 0pt \leftskip 0pt \rightskip 0pt plus 1fil \def\\{\par}
\sl\theaddress\par
\medskip
\rm Email:\stdspace\tt\theemail\hfill\rm Received:\qua\receiveddate \par}}

\def\recd{{\small Received:\qua\receiveddate\ifx\reviseddate\relax
\else\qquad Revised:\qua\reviseddate\fi\par}} 

\def\lognumber#1{\def\thelognumber{#1}}
\def\volumenumber#1{\def\thevolumenumber{#1}}
\def\volumeyear#1{\def\thevolumeyear{#1}}
\def\papernumber#1{\def\thepapernumber{#1}}
\def\pagenumbers#1#2{\def\startpage{#1}\def\finishpage{#2}}
\def\published#1{\def\publishdate{#1}}

\def\received#1{\def\receiveddate{#1}}

\def\accepted#1{\def\accepteddate{#1}}

\long\def\asciiabstract#1{\long\def\theasciiabstract{#1}}

\let\\\par\let\thelognumber\relax\let\thevolumenumber\relax
\let\thepapernumber\relax\let\thevolumeyear\relax\let\startpage\relax
\let\finishpage\relax\let\publishdate\relax\let\receiveddate\relax
\let\reviseddate\relax\let\accepteddate\relax\let\theasciititle\relax
\let\theasciiauthors\relax
\let\theasciiabstract\relax

\let\theasciiemail\relax

\ifplaintex
\font\logobig=cmssbx10 scaled 3836
\font\logomed=cmssbx10 scaled 2557
\else
\font\logobig=cmssbx10 scaled 4200
\font\logomed=cmssbx10 scaled 2800
\fi

\long\def\makeagttitle{   
\count0=\startpage
\agt\hfill      
\hbox to 45truept{\vbox to 0pt{\vglue -13truept{\logomed A\kern -.37em{\logobig 
T}\kern -.38em G}\vss}\hss}
\break
{\small Volume \thevolumenumber\ (\thevolumeyear)
\startpage--\finishpage\nl
Published: \publishdate}

\vglue .25truein

{\parskip=0pt\leftskip 0pt plus
1fil\def\\{\par\smallskip}{\Large\bf\thetitle}\par\medskip} \vglue
0.05truein

{\parskip=0pt\leftskip 0pt plus 1fil\def\\{\par}{\sc\theauthors}
\par\medskip}%
 
\vglue 0.03truein 

{\small\leftskip 25truept\rightskip 25truept{\bf Abstract}\stdspace\theabstract

{\bf AMS Classification}\stdspace\theprimaryclass
\ifx\thesecondaryclass\relax\else; \thesecondaryclass\fi\par
{\bf Keywords}\stdspace \thekeywords\par}\vglue 7truept

}   

\ifplaintex
\font\phead=cmsl9 scaled 950
\font\pnum=cmbx10 scaled 913
\font\pfoot=cmsl9 scaled 950
\headline{\vbox to 0pt{\vskip -4.5mm\line{\small\phead\ifnum
\count0=\startpage ISSN 1472-2739 (on-line) 1472-2747 (printed)
\hfill {\pnum\folio}\else\ifodd\count0\def\\{ }%
\ifx\theshorttitle\relax\thetitle\else\theshorttitle\fi\hfill{\pnum\folio}
\else\def\\{ and }{\pnum\folio}\hfill\ifx\theshortauthors\relax\theauthors
\else\theshortauthors\fi\fi\fi}\vss}}
\footline{\vbox to 0pt{\vglue 0mm\line{\small\pfoot\ifnum\count0=\startpage
\copyright\ \gtp\hfill\else
\agt, Volume \thevolumenumber\ (\thevolumeyear)\hfill\fi}\vss}}
\else
\font=cmsl9 scaled 1050
\font\lnum=cmbx10 
\font=cmsl9 scaled 1050
\makeatletter
\def\@oddhead{{\small\lhead\ifnum\count0=\startpage ISSN 1472-2739 
(on-line) 1472-2747 (printed)\hfill {\lnum\number\count0}\else\ifodd\count0
\def\\{ }\ifx\theshorttitle\relax \thetitle \else\theshorttitle\fi\hfill
{\lnum\number\count0}\else\def\\{ and }{\lnum\number\count0}
\hfill\ifx\theshortauthors\relax 
\theauthors\else\theshortauthors\fi\fi\fi}}\def\@evenhead{\@oddhead}
\def\@oddfoot{\small\lfoot\ifnum\count0=\startpage\copyright\ \gtp\hfill\else
\agt, Volume \thevolumenumber\ (\thevolumeyear)\hfill\fi}
\def\@evenfoot{\@oddfoot}
\makeatother
\fi
\let\maketitlepage\makeagttitle
\let\makeshorttitle\maketitlepage
\let\maketitle\maketitlepage

\newwrite\gtoutfile
\long\gdef\makeheadfile{  
{\def\\{, }\def\s{ }
\immediate\openout\gtoutfile head.xxx
\immediate\write\gtoutfile{To: math@arxiv.org}
\immediate\write\gtoutfile{Subject: put OR rep NNNNN:ppppp}
\immediate\write\gtoutfile{--text follows this line--}
\immediate\write\gtoutfile{Proxy-for: \ifx\theasciiauthors\relax
\theauthors\else\theasciiauthors\fi\s<\ifx\theasciiemail\relax\theemail\else\theasciiemail\fi>}
\immediate\write\gtoutfile{\noexpand\\}
\immediate\write\gtoutfile{Authors: \ifx\theasciiauthors\relax
\theauthors\else\theasciiauthors\fi}
{\def\\{ }\immediate\write\gtoutfile{Title: \ifx\theasciititle\relax
\thetitle\else\theasciititle\fi}}
\immediate\write\gtoutfile{Subj-class: GT or SG, GR etc}
\immediate\write\gtoutfile{MSC-class: \theprimaryclass\ifx\thesecondaryclass\relax\else, \thesecondaryclass\fi}
\immediate\write\gtoutfile{Journal-ref: Algebr. Geom. Topol. \thevolumenumber\s
(\thevolumeyear) \startpage-\finishpage}
\immediate\write\gtoutfile{Comments: Published by Algebraic and
Geometric Topology at}
\immediate\write\gtoutfile{\s\s\s  http://www.maths.warwick.ac.uk/agt/AGTVol\thevolumenumber/agt-\thevolumenumber-\thepapernumber.abs.html}
\immediate\write\gtoutfile{\noexpand\\}
\immediate\write\gtoutfile{}
\ifx\theasciiabstract\relax
\immediate\write\gtoutfile{\theabstract}\else
\immediate\write\gtoutfile{\theasciiabstract}\fi
\immediate\write\gtoutfile{}
\immediate\write\gtoutfile{\noexpand\\}
\immediate\write\gtoutfile{}
\immediate\closeout\gtoutfile}}  

\def\maketitlepage{\makeagttitle\makeheadfile}
\let\makeshorttitle\maketitlepage
\let\maketitle\maketitlepage

\def\ifplaintex{\expandafter\ifx\csname documentclass\endcsname\relax}

\def\gtp{{\mathsurround=0pt\it $\cal G\mskip-2mu$eometry \&\ 
$\cal T\!\!$opology $\cal P\!$ublications}}  

\def\Addressesr{\bigskip
{\small \parskip 0pt \leftskip 0pt \rightskip 0pt plus 1fil \def\\{\par}
\sl\theaddress\par
\medskip
\rm Email:\stdspace\tt\theemail\hfill\rm Received:\qua\receiveddate \par}}

\def\recd{{\small Received:\qua\receiveddate\ifx\reviseddate\relax
\else\qquad Revised:\qua\reviseddate\fi\par}} 

\def\lognumber#1{\def\thelognumber{#1}}
\def\volumenumber#1{\def\thevolumenumber{#1}}
\def\volumeyear#1{\def\thevolumeyear{#1}}
\def\papernumber#1{\def\thepapernumber{#1}}
\def\pagenumbers#1#2{\def\startpage{#1}\def\finishpage{#2}}
\def\published#1{\def\publishdate{#1}}

\def\received#1{\def\receiveddate{#1}}

\def\accepted#1{\def\accepteddate{#1}}

\long\def\asciiabstract#1{\long\def\theasciiabstract{#1}}

\let\\\par\let\thelognumber\relax\let\thevolumenumber\relax
\let\thepapernumber\relax\let\thevolumeyear\relax\let\startpage\relax
\let\finishpage\relax\let\publishdate\relax\let\receiveddate\relax
\let\reviseddate\relax\let\accepteddate\relax\let\theasciititle\relax
\let\theasciiauthors\relax
\let\theasciiabstract\relax

\let\theasciiemail\relax

\ifplaintex
\font\logobig=cmssbx10 scaled 3836
\font\logomed=cmssbx10 scaled 2557
\else
\font\logobig=cmssbx10 scaled 4200
\font\logomed=cmssbx10 scaled 2800
\fi

\long\def\makeagttitle{   
\count0=\startpage
\agt\hfill      
\hbox to 45truept{\vbox to 0pt{\vglue -13truept{\logomed A\kern -.37em{\logobig 
T}\kern -.38em G}\vss}\hss}
\break
{\small Volume \thevolumenumber\ (\thevolumeyear)
\startpage--\finishpage\nl
Published: \publishdate}

\vglue .25truein

{\parskip=0pt\leftskip 0pt plus
1fil\def\\{\par\smallskip}{\Large\bf\thetitle}\par\medskip} \vglue
0.05truein

{\parskip=0pt\leftskip 0pt plus 1fil\def\\{\par}{\sc\theauthors}
\par\medskip}%
 
\vglue 0.03truein 

{\small\leftskip 25truept\rightskip 25truept{\bf Abstract}\stdspace\theabstract

{\bf AMS Classification}\stdspace\theprimaryclass
\ifx\thesecondaryclass\relax\else; \thesecondaryclass\fi\par
{\bf Keywords}\stdspace \thekeywords\par}\vglue 7truept

}   

\ifplaintex
\font\phead=cmsl9 scaled 950
\font\pnum=cmbx10 scaled 913
\font\pfoot=cmsl9 scaled 950
\headline{\vbox to 0pt{\vskip -4.5mm\line{\small\phead\ifnum
\count0=\startpage ISSN 1472-2739 (on-line) 1472-2747 (printed)
\hfill {\pnum\folio}\else\ifodd\count0\def\\{ }%
\ifx\theshorttitle\relax\thetitle\else\theshorttitle\fi\hfill{\pnum\folio}
\else\def\\{ and }{\pnum\folio}\hfill\ifx\theshortauthors\relax\theauthors
\else\theshortauthors\fi\fi\fi}\vss}}
\footline{\vbox to 0pt{\vglue 0mm\line{\small\pfoot\ifnum\count0=\startpage
\copyright\ \gtp\hfill\else
\agt, Volume \thevolumenumber\ (\thevolumeyear)\hfill\fi}\vss}}
\else
\font=cmsl9 scaled 1050
\font\lnum=cmbx10 
\font=cmsl9 scaled 1050
\makeatletter
\def\@oddhead{{\small\lhead\ifnum\count0=\startpage ISSN 1472-2739 
(on-line) 1472-2747 (printed)\hfill {\lnum\number\count0}\else\ifodd\count0
\def\\{ }\ifx\theshorttitle\relax \thetitle \else\theshorttitle\fi\hfill
{\lnum\number\count0}\else\def\\{ and }{\lnum\number\count0}
\hfill\ifx\theshortauthors\relax 
\theauthors\else\theshortauthors\fi\fi\fi}}\def\@evenhead{\@oddhead}
\def\@oddfoot{\small\lfoot\ifnum\count0=\startpage\copyright\ \gtp\hfill\else
\agt, Volume \thevolumenumber\ (\thevolumeyear)\hfill\fi}
\def\@evenfoot{\@oddfoot}
\makeatother
\fi
\let\maketitlepage\makeagttitle
\let\makeshorttitle\maketitlepage
\let\maketitle\maketitlepage

\newwrite\gtoutfile
\long\gdef\makeheadfile{  
{\def\\{, }\def\s{ }
\immediate\openout\gtoutfile head.xxx
\immediate\write\gtoutfile{To: math@arxiv.org}
\immediate\write\gtoutfile{Subject: put OR rep NNNNN:ppppp}
\immediate\write\gtoutfile{--text follows this line--}
\immediate\write\gtoutfile{Proxy-for: \ifx\theasciiauthors\relax
\theauthors\else\theasciiauthors\fi\s<\ifx\theasciiemail\relax\theemail\else\theasciiemail\fi>}
\immediate\write\gtoutfile{\noexpand\\}
\immediate\write\gtoutfile{Authors: \ifx\theasciiauthors\relax
\theauthors\else\theasciiauthors\fi}
{\def\\{ }\immediate\write\gtoutfile{Title: \ifx\theasciititle\relax
\thetitle\else\theasciititle\fi}}
\immediate\write\gtoutfile{Subj-class: GT or SG, GR etc}
\immediate\write\gtoutfile{MSC-class: \theprimaryclass\ifx\thesecondaryclass\relax\else, \thesecondaryclass\fi}
\immediate\write\gtoutfile{Journal-ref: Algebr. Geom. Topol. \thevolumenumber\s
(\thevolumeyear) \startpage-\finishpage}
\immediate\write\gtoutfile{Comments: Published by Algebraic and
Geometric Topology at}
\immediate\write\gtoutfile{\s\s\s  http://www.maths.warwick.ac.uk/agt/AGTVol\thevolumenumber/agt-\thevolumenumber-\thepapernumber.abs.html}
\immediate\write\gtoutfile{\noexpand\\}
\immediate\write\gtoutfile{}
\ifx\theasciiabstract\relax
\immediate\write\gtoutfile{\theabstract}\else
\immediate\write\gtoutfile{\theasciiabstract}\fi
\immediate\write\gtoutfile{}
\immediate\write\gtoutfile{\noexpand\\}
\immediate\write\gtoutfile{}
\immediate\closeout\gtoutfile}}  

\def\maketitlepage{\makeagttitle\makeheadfile}
\let\makeshorttitle\maketitlepage
\let\maketitle\maketitlepage

\lognumber{41}
\volumenumber{2}
\volumeyear{2002}
\papernumber{41}
\pagenumbers{1051}{1060}
\received{13 June 2002}
\accepted{29 Oct 2002}
\published{5 November 2002}

\usepackage{amsmath, epsf}  

\def\zz{{\bf Z}}
\def\qq{{\bf Q}}

\def\<{\langle}
\def\>{\rangle}

\newtheorem{theorem}{Theorem}[section]
\newtheorem{lemma}[theorem]{Lemma}

\newtheorem{conjecture}[theorem]{Conjecture}

\theoremstyle{definition}
\newtheorem{definition}[theorem]{Definition}

\theoremstyle{remark}

\numberwithin{equation}{section}

\begin{document}
 
\title{The slicing number of a knot}

\author{Charles Livingston}

\address{Department of Mathematics, Indiana University, Bloomington, IN 47405, USA}

\email{livingst@indiana.edu}

\begin{abstract}  An open question asks if every knot of 4--genus
$g_s$ can be changed into a slice knot by $g_s$ crossing changes.  A
counterexample is given. 
\end{abstract}

\asciiabstract{An open question asks if every knot of 4-genus
g_s can be changed into a slice knot by g_s crossing changes.  A
counterexample is given.}



\maketitle

A question of Askitas, appearing  in \cite[Problem 12.1]{o}, asks the
following: Can a knot of 4--genus
$g_s$ always be sliced (made into a slice knot) by $g_s$ crossing changes?  If
we let
$u_s(K)$ denote the {\em slicing number} of $K$, that is, the minimum number of
crossing changes that are needed to convert
$K$ into a slice knot, one readily shows that $g_s(K) \le u_s(K)$ for all
knots, with equality if $g_s(K) = 0$.  Hence, the problem can be restated as
asking if  $g_s(K) = u_s(K)$ for all
$K$.

We will show that the knot
$7_4$ provides a counterexample; $g_s(7_4) = 1$ but no
crossing change results in a slice knot: $u_s(7_4) = 2$. It is interesting to
note that     $7_4$ already stands out as an important example. The proof
that its unknotting number is 2, not 1,   resisted early attempts 
\cite{n};   ultimately,   Lickorish   \cite{l} succeeded in proving that it
cannot be unknotted with a single crossing change. 

 As noted by Stoimenow in \cite{o}, if one attempts to unknot a knot of 4--ball
genus $g_s$ instead of converting it into a slice knot, more than $g_s$
crossing changes may be required. This is obviously the case with slice knots. 
For a more general example, let   $T$ denote the
 trefoil knot. One has $g_s( n(T \# -T) \# mT) = m$, but $u( n(T
\# -T) \# mT) = m + 2n$, where $u$ denotes the unknotting number.  

The unknotting number, though itself mysterious, appears much simpler than the
slicing number.  Many of the three--dimensional tools that are available for
studying the unknotting number do not apply to the study of the slicing number.
As we will see,   even for this low crossing knot,
$7_4$, the computation of its slicing number is far more complicated than its
unknotting number.

In the last section of this paper   we introduce a new slicing invariant,
$U_s(K)$, that takes into account the sign of crossing changes used to convert
a knot $K$ into a slice knot.  This invariant is more closely related to the
4--genus and satisfies
 $$g_s(K)
\le U_s(K) \le u_s(K).$$  It seems  likely that there are knots $K$ for which
$g_s(K) \ne U_s(K)$, and $7_4$ seems a good candidate, but we have been
unable to verify this.

A good  reference  for the knot theory used here, especially surgery
descriptions of knots, crossing changes and branched coverings,   is
\cite{r}.  A reference for   4--dimensional aspects of knotting and also for
the linking form of 3--manifolds is
\cite{go}.  A careful analysis of the interplay between crossing changes
and the linking form of the 2--fold branched cover of a knot appears in  
\cite{l},   which our work here generalizes.  Different
aspects of the relationship between crossing changes and 4--dimensional aspects
of knotting appear in \cite{cl}.  A general discussion of slicing operations is
contained in \cite{a}.

\section{Background} Our goal  is to prove that a single crossing change cannot
change $7_4$ into a slice knot.  The key results  concerning slice knots that
we will be using are contained in  the following theorem; details of the proof
can be found in \cite{cg, go, r}.

\begin{theorem} \label{condition} If $K$ is slice then: 

\begin{enumerate}

\item   $\Delta_K(t) =
\pm f(t)f(t^{-1})$ for some polynomial $f$, where $\Delta_K(t) $ is the
Alexander polynomial.

\item $|H_1(M(K),\zz)| = n^2$ for some odd $n$, where $M(K)$ is the 2--fold
branched cover of $S^3$ branched over $K$.

\item There is a subgroup $H \subset H_1(M(K),\zz)$ such that $|H|^2 =
|H_1(M(K),\zz)|$ and the $\qq/\zz$--valued linking form $\beta$ defined on $
H_1(M(K),\zz)$ vanishes on $H$.
\end{enumerate}
\end{theorem}

 Our analysis of $7_4$ will focus on the 2--fold branched cover, $M(7_4)$, and
its linking form.  This is much as in Lickorish's unknotting number argument. 
However, in our case the necessary  analysis of the 2--fold branched cover can
only be achieved by a close examination of the infinite cyclic cover. 
In the next two sections we examine the 2--fold branched cover; in Section
\ref{infinite}   we consider the infinite cyclic cover.


\section{Crossing Changes and Surgery}
\label{crossings}

If a knot $K'$ is obtained from $K$ by changing a crossing,   surgery theory as
described in 
\cite{r} quickly gives that the 2-fold cover of $K'$, $M(K')$, can be obtained
from $M(K)$ by performing integral surgery on a pair curves, say $S_1$ and
$S_2$, in $S^3$.  It is  also known \cite{l, mo} that
$M(K')$ can be obtained from $M(K)$ by performing $p/2$ surgery on a single
curve, say $T$, in $S^3$.  Here it will be useful to observe that $T$ can be
taken to be
$S_1$, as we next describe.

A crossing change is formally achieved as follows.  Let $D$ be a disk meeting
$K$ transversely in two points.  A neighborhood of $D$ is homeomorphic to a
3--ball, $B$, meeting $K$ in two trivial arcs. In one view, a crossing change
is accomplished by performing
$\pm 1$ surgery on the boundary of $D$, say $S$.  Then $S$ lifts to give the
curves $S_1$ and $S_2$ in $M(K)$.  In the other view, the crossing change is
accomplished by removing $B$ from $S^3$ and sewing it back in with one full 
twist.  The 2--fold branched cover of $B$ is a single solid torus, a regular
neighborhood of its core $T$. A close examination shows the surgery coefficient
in this case is $p/2$ for some odd $p$. 

The lift of $D$ to the 2-fold branched cover is an annulus with boundary the
union of $S_1$ and $S_2$ and core $T$ (the lift of an arc, $\tau$, on $D$ with
endpoints the two points of intersection of $D$ with $K$).  Clearly $T$ is
isotopic to either $S_i$, as desired.

The following generalization of these observations will be useful.  Rather than
put a single full twist between the strands when replacing $B$, $n$ full twists
can be added.  This is achieved by performing $\pm 1/n$ surgery on $S$ and hence
the 2--fold branched cover is obtained by performing $p/n$ surgery on the $S_i$
for some $p$, or, by a similar analysis, by performing $p'/2n$ surgery on $T$
for some $p'$.


\section{Results based on the 2--fold branched cover of  $7_4$ }\label{2fold}

\begin{figure}[ht!]
$$ \epsfxsize=1.2in  \epsfbox{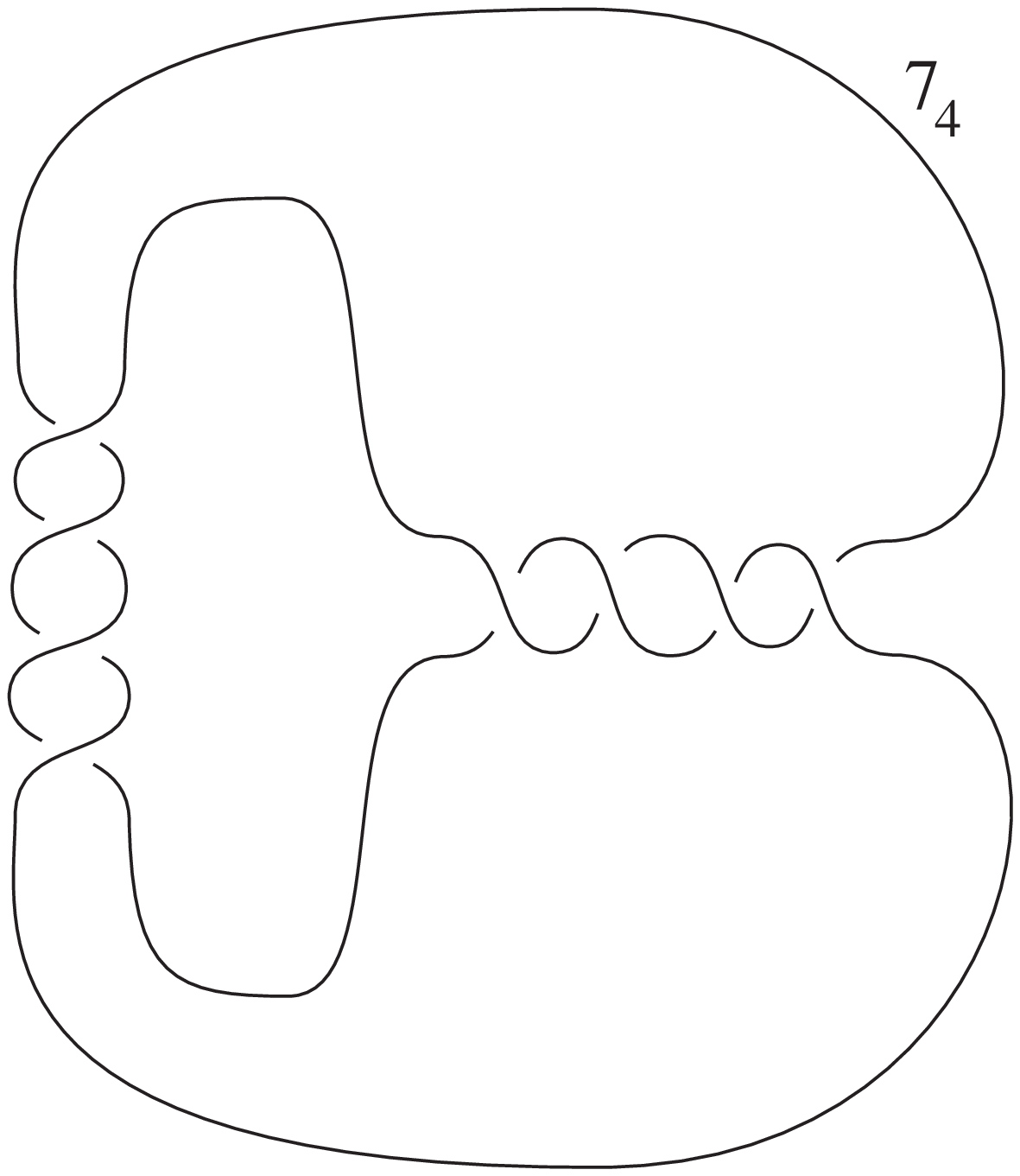} \hskip.9in  \epsfxsize=1.9in 
\epsfbox{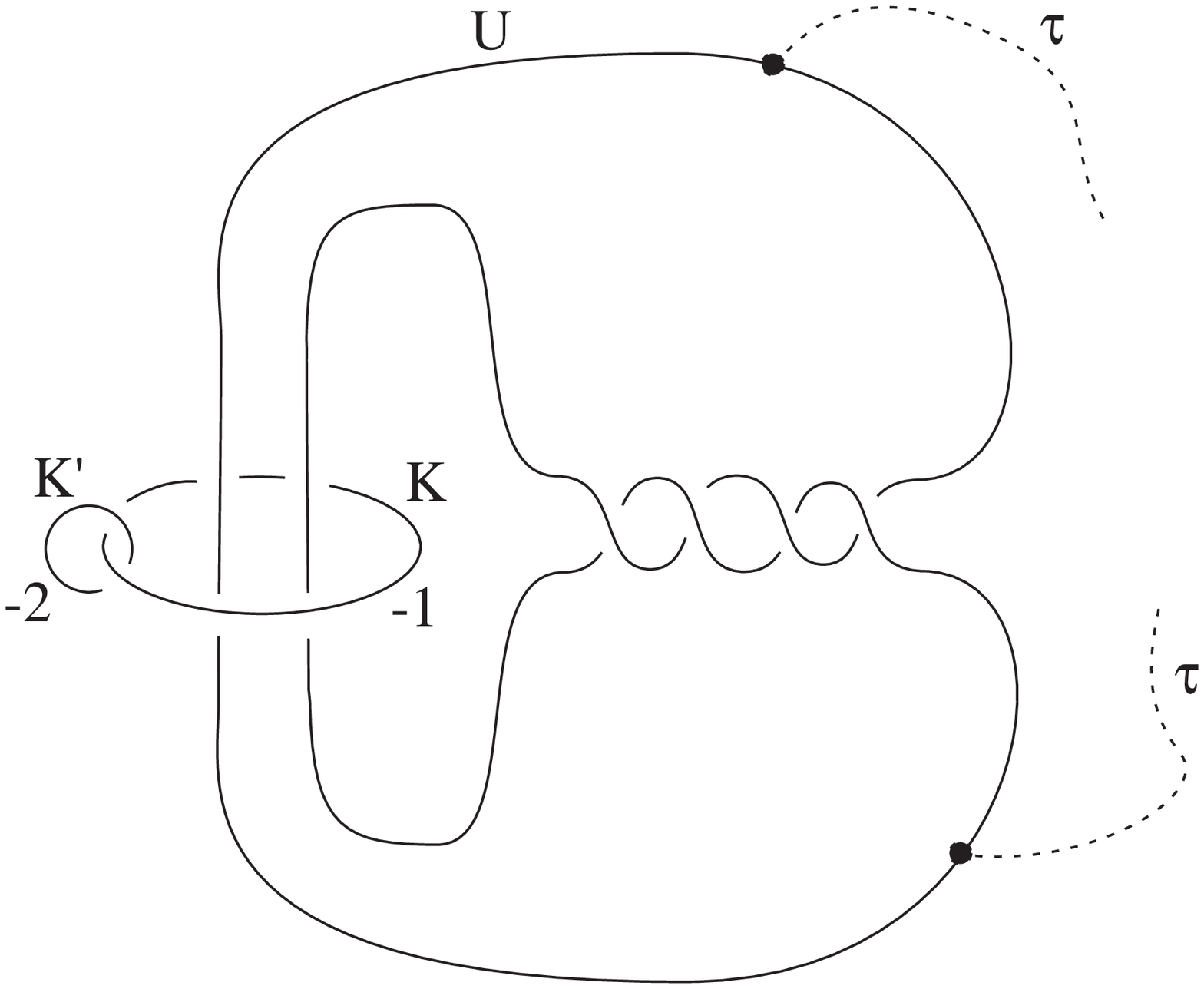} 
$$
\caption{The knot $7_4$} \label{7_4figure}
\end{figure}

On the left in  Figure \ref{7_4figure} the knot $7_4$ is illustrated.   Basic
facts about $7_4$   include that it has 3--sphere genus 1 and that its
Alexander polynomial is $\Delta_{7_4} = 4t^2 - 7t +4$.  Since the Alexander
polynomial is irreducible,  
$7_4$ is not slice, so we have $g_s(7_4) = 1$.   Also,
$7_4$ is the 2--bridge   $B(4,-4)$, and hence from the continued fraction
expansion  it has 2--fold branched cover the lens space, $L(15,4)$.

The right diagram in Figure \ref{7_4figure} represents a surgery diagram of
$7_4$.  According to
\cite{r}, surgery on the link $K \amalg K'$ with coefficient $-1$ and $-2$
yields $S^3$. Also according to \cite{r} the component $K'$ could be ignored in
the diagram if $-1/2$ surgery is performed on $K$ instead.  In both cases the
effect is to put two full right handed twists in the two strands passing
through $K$.  

Notice  that 
$U$ is unknotted.   After surgery is performed, $U$ is converted into the knot
$7_4$.

 If a knot $J$ is obtained from $7_4$ by a single crossing change, that  change
is achieved via a disk $D$ meeting $7_4$ in two points, marked schematically by
the two dots in the right hand diagram.  The path on
$D$ joining those two points is denoted
$\tau$, a portion of which is also indicated schematically. By sliding $\tau$
over $K$ repeatedly it can be arranged that $\tau$ misses the small disk
bounded by $K'$ meeting $K$ in one point. The boundary of
$D$ will be denoted $S$ and one of its lift to the 2--fold branched cover of
$S^3$ over $U$ (this cover is again
$S^3$ since $U$ is unknotted) will be   denoted $S_1$.  Neither $D$ nor $S$ is
drawn in the figure.

Since two full twists on the unknot $U$ convert it into $7_4$, the 2--fold
branched cover of $S^3$ branched over
$7_4$ is, by our earlier discussion, obtained from $S^3$ by surgery on a single
lift of $K$, say
$K_1$, with surgery coefficient  of the form $p/4$ for some $p$.  Since we know
that the cover is
$L(15,4)$, we actually know that $p = 15 $, though for the argument that
follows, simply knowing that
$p = \pm 15$ would be sufficient.

\begin{theorem}\label{theorem1} If the linking number of $K_1$ and $S_1$ in
$S^3$ is divisible by 15 then
$J$ is not slice.
\end{theorem}

\begin{proof} Suppose that the linking number is divisible by 15.  Since $15/4$
surgery is performed on
$K_1$, after repeatedly sliding
$S_1$ over $K_1$  it can be arranged that the linking number of $K_1$ and $S_1$
is 0.  The 2--fold cover of $S^3$ branched over $J$, that is  $M(J)$, is
obtained from $S^3$ by performing
$  15/4$ surgery on
$K_1$ and $p/2$ surgery on $S_1$ for some odd $p$.

If $J$ is slice, the order of the homology of  $M(J)$ is an odd square and
hence $ p =
\pm  5^{2k+1} 3^{2j+1} q^2$, where $q$ is relatively prime to 30.

We have that $H_1(M(J),\zz) = \zz_{15} \oplus \zz_{|p|}$ generated by the
meridians of $K_1$ and $S_1$, denoted $m_1 $ and $m_2$, respectively.  

The $\qq/\zz$--valued linking form, $\beta$, on $H_1(M(J),\zz)$ is orthogonal
with respect to this direct sum decomposition since the linking number is now
0.  Furthermore, from the surgery description we have that $\beta( m_1, m_1) =
  4/15$ and $\beta(m_2,m_2) = 2/p$.  The 5--torsion in $H_1(M(J),\zz)$ is
isomorphic to $\zz_5 \oplus
\zz_{5^{2k+1}}$, generated by $n_1 = 3 m_1$ and $n_2 = 3^{2j+1}q^2m_2$.  A
quick calculation shows that $\beta(n_1,n_1) =   2/5$ and $\beta(m_2,m_2) = 2
(3^{2j+1} q^2)^2/p = 
\pm 2 (3^{2j+1} q^2) /5^{2k+1}$.

If $J$ is slice, the linking form on the 5--torsion vanishes on a subgroup of
order
$5^{k+1}$.   Suppose that $n_1 + x5^l n_2$ has self--linking
$0 \in \qq/\zz$, where $x$ is relatively prime to $5$.  Then we would have
$$ \frac{2}{ 5}
\pm 
\frac{2 x^2 5^{2l}(3^{2j+1} q^2)}{5^{2k+1} }= 0 \in \qq/\zz.$$ This implies
that $l = k$, and hence that $  2 \pm 2 x^23^{2j+1} q^2 \equiv 0
\ \mbox{mod}\  5$.  Letting $q' = x3^jq$, this can be rewritten as $\  2 \pm
2(3{q'}^2) \equiv 0
\ \mbox{mod}\  5$, or that $2 \equiv \mp {q'}^2 \mbox{mod}\ 5$.  However, the
only squares modulo 5 are
$\pm 1$, so this is impossible.

It follows from this that any element of self--linking 0 must be of the form
$x5^l n_2$ for some $l$ and $x$ relatively prime to 5.  One quickly computes
that $l > k$, but such elements  generate a subgroup of order $5^{k}$, which is
not large enough to satisfy the condition of Theorem
\ref{condition}, Statement 3.
\end{proof}


\section{The Infinite Cyclic Cover of $7_4$ }\label{infinite}
 The goal of this section is to prove the following result.  It, along with
Theorem
\ref{theorem1}, shows that $7_4$ cannot be changed into a slice knot with a
single crossing change.
\begin{theorem}\label{theorem2}  If a crossing change converts $7_4$ into a
slice knot $J$, then the corresponding curve $S_1$ in $M(7_4)$ is null
homologous in $H_1(L(15,4),\zz)$.
\end{theorem}

Before beginning the proof   we need to set up notation and prove
a lemma.

The infinite cyclic cover of $J$ is built from the infinite cyclic cover of the
unknot, $U$, by performing equivariant surgery on three families of curves:
$\{\tilde{K}_i\},  \{\tilde{K_i'}\}$ and $\{\tilde{S}_i\}$, using the notation
as before.  (In each case, $i = -\infty, \dots , \infty$.)

Following Rolfsen \cite{r}, one can draw that cover with the $\{\tilde{K}_i\}, 
\{\tilde{K_i'}\}$ drawn explicitly, and the $\{\tilde{S}_i\}$ unknown curves.
From this one finds the presentation matrix of the infinite cyclic cover of $J$
as a $\zz[t,t^{-1}]$ module, with respect to the basis given by the meridians
of  $ \tilde{K}_0 , 
 \tilde{K_0'}$ and $ \tilde{S}_0$, say $k_0, k'_0$, and $s_0$. The resulting
presentation is given by the matrix
$$ A = \left( \begin{matrix} -2t + 3 - 2t^{-1} & 1 & g(t)\\
               1 & -2 & 0 \\
               g(t^{-1}) &  0 &  f(t)
\end{matrix} \right). 
$$ Here $g(t)$ is an unknown polynomial describing the linking between the
lifts of $S$ and those of $K$.  (Notice that the lifts of $S$ do not link the
lifts of $K'$, since $\tau$ (and so
$S$) misses the small disk bounded by $K'$ and this disk lifts to a series of
disjoint disks bounded by the $\tilde{K_i'}$ in the infinite cyclic cover.) 
Also,
$f(t)$ is an unknown symmetric polynomial describing the self--linking of the
lifts of $S$. (It might be helpful for the reader to note that if $g = 0$ and
$f = 1$ then the determinant of the  matrix is $4t - 7 + 4t^{-1}$, the
Alexander polynomial of $7_4$.)

Although $g$ and $f$ are  unknown, two observations are possible.  The first is
that $f(1) =
\pm 1$; this is because   $\pm 1$ surgery is being performed on $S$.  The
second is that
$g(1) = 0$, or that $(t-1) $ divides $g$, which follows from the fact that $S$
and $K$ have 0 linking number, since $S$ bounds the disk $D$ in the complement
of $K$.

\begin{lemma} If $J$ is slice, then  $4t- 7 +4t^{-1} = \Delta_{7_4}$  divides
$g$. 

\end{lemma}

\begin{proof} The determinant of
$A$ is given by 
$$\Delta_J(t) = f(t)\Delta_{7_4}(t) + 2g(t)g(t^{-1}).$$  Since $J$ is assumed
to be slice we can rewrite this as $$\pm H(t)H(t^{-1}) = f(t)\Delta_{7_4}(t) +
2g(t)g(t^{-1})$$ for some
$H(t)$.  Clearly, if $H$ is divisible by $\Delta_{7_4}$ then $g(t)$ would also
be and we would be done.  So, assume that neither $H$ or $g$ has factor
$\Delta_{7_4}$. 

Working modulo $\Delta_{7_4}$ we now have the equation:
$$(*)\hskip.4in 2g(t)g(t^{-1}) = \pm H(t)H(t^{-1})\ \in \ \zz[t,t^{-1}]/\<4t - 7
+ 4t^{-1}\>.$$
 There is an injection $\phi: \zz[t,t^{-1}]/\<4t - 7 + 4t^{-1}\> \to
\qq(\sqrt{-15})$ with
$\phi(t^{\pm 1} ) = (7 \pm \sqrt{-15})/4$. It follows that if equation $(*)$
holds then we could factor $2 = \pm (\frac{a}{c} + \frac{b}{c}\sqrt{-15} )( 
(\frac{a}{c} -
\frac{b}{c}\sqrt{-15} )$ with $a, b$, and $c$ integers with $\gcd(a,b,c) = 1$.
Simplifying we would have $$\pm 2c^2 - a^2 - 15b^2 = 0.$$
Working modulo 5 and using that $\pm 2$ is not a quadratic residue modulo 5, one
sees immediately that
$a$ and
$c$ are both divisible by 5, which implies (working modulo 25) that $b$ is
divisible by 5 as well.  Write
$a = 5^s a'$,  $b = 5^t b'$ and $c = 5^r c'$, with $a', b',$ and $c'$
relatively prime to 5.  Hence:
$$\pm 2( 5^{2s}{c'}^2) - 5^{2t}{a'}^2 - 3( 5^{2r+1}{b'}^2 )= 0.$$
If among the three exponents of 5 that appear in this equation there is a
unique smallest exponent, then factoring out that power of 5 leaves an equation
that clearly cannot hold modulo 5.  Hence, there must be two exponents that are
equal, and these must be the two even exponents.  Factoring these out leaves
the equation:
$$\pm 2  {c'}^2  -  {a'}^2 - 3( 5^{2r'+1}{b'}^2 )= 0.$$  Again using that $\pm
2$ is not a quadratic residue modulo 5 gives a contradiction. 
\end{proof}

We can now prove Theorem \ref{theorem2}.

\begin{proof}[Proof of  Theorem \ref{theorem2}] The polynomial $g$ determines
the linking numbers of the lifts of $K$ and $S$ to the $n$-fold cyclic branched
cover of $S^3$ branched over $U$ as follows.  Call the  lifts $\bar{K}_i$ and
$\bar{S}_i$ with $i$ running from 0 to
$n-1$. The linking numbers are given by equivariance and 
$$\mbox{lk}( \bar{K}_0 , \bar{S}_i) = \bar{g}_i$$ where $\bar{g}_i$ is the
coefficient of
$t^i$ in the reduction $\bar{g}$ of $g$ to $\zz[t,t^{-1}] /  \<t^n - 1\>$.

In the case of the 2--fold cover we are hence interested in the even and odd
index coefficients.  For any integral polynomial
$F(x) =
\sum a_i t^i$ the sum of the even index coefficients is given by
$(F(1) + F(-1))/2$ and the sum of the odd index coefficients is $(F(1) -
F(-1))/2$.  In our case we have seen that $g(t) = (t-1)(4t^2 - 7t + 4) h(t) $
for some $h$.  Hence, the sum of the even (or odd) coefficients is given by
$\pm 15 h(-1)$.  In particular, the linking number is divisible by 15.  Hence
$\bar{S}_i = S_i$ is null homologous in the $L(15,4)$ obtained by surgery on
$K_1$.
\end{proof}


\section{Extensions}\label{extend}

The proof that $7_4$ has slicing number 2 clearly generalizes to other knots,
though a general statement is somewhat technical.  On the other hand, these
methods seem not to apply effectively in addressing the next level of
complexity---finding a knot $K$ with $g_s(K) = 2$ but with slicing number 3.  

\begin{conjecture}The difference $u_s(K) - g_s(K)$ can be arbitrarily large.  
\end{conjecture} \noindent In fact, this gap should be arbitrarily large even
for knots with $g_s = 1$.

In retrospect, Askitas's question was optimistic.  It is easily seen that if a
knot can be converted into a slice knot by making $n$ positive and $n$ negative
crossing changes, then
$g_s(K) \le n$. More generally, we have the following signed unknotting number.

\begin{definition} For a knot $K$, let $I$ denote the set of pairs of
nonnegative integers $(m,n)$ such that some collection of $m$ positive crossing
changes and $n$ negative crossing changes  
  converts $K$ into a slice knot. Define the invariant $U_s(K)$ by $$U_s(K) =
\min_{(m,n) \in I}\{\max(m,n)\}$$
\end{definition}

The following result has an elementary proof.
 \begin{theorem} For all $K$, $g_s(K) \le U_s(K)$.\end{theorem}

The only bounds that we know of relating to $U_s$ are those arising from $g_s$,
and so it is possible that $U_s(K) = g_s(K)$ for all $K$.  However, a more
likely conjecture is the following.

\begin{conjecture} The difference $  U_s(K) - g_s(K) $ can be arbitrarily large.
\end{conjecture}

Even the following example is unknown.
\vskip.1in
\noindent{\bf Question}\qua Does  $U_s(7_4) = 1$?
\vskip.1in
The example we describe below indicates that proving that $U_s(7_4) = 2$ may be
quite difficult.

\vskip.2in

\noindent {\bf General Twisting}\qua  One can think of performing a crossing change
as grabbing two parallel strands of a knot with opposite orientation  and given
them one full twist.  More generally, one can grab $2k$ parallel strands of $K$
with $k$ of the strands oriented in each direction and giving them one full
twist.   Call this a {\em generalized crossing change}.  With a little care,
the proof that
$7_4$ cannot be converted into a slice knot generalizes to show the  following:

\begin{theorem} The knot $7_4$ cannot be converted into a slice knot using a
single generalized crossing change. \end{theorem}

\begin{figure}[ht!]
\cl{\epsfxsize=2in  \epsfbox{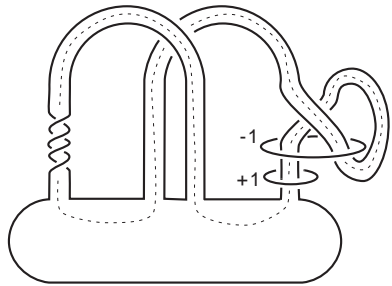}}
\caption{Twisting $7_4$ to a slice knot} \label{twist74}
\end{figure}

On the other hand, consider Figure \ref{twist74}. The illustrated knot is slice
since the dotted curve on the Seifert surface is unknotted and has framing 0. 
If a right--handed twist is put on the strands going through the circle
labelled $-1$ and a left--handed twist is put on the strands going through the
circle labelled $+1$, then the knot $7_4$ results.  Hence, $7_4$ can be
converted into a slice knot by performing one positive and one negative {\em
generalized} crossing change.  

Since all the relevant techniques that we know of do not distinguish between
crossing changes and generalized crossing changes, the difficulty associated to
disproving showing that $U_s(7_4) = 2$ is now clear.

It is worth pointing out here that clearly $7_4$ can be converted into a slice
knot (actually the unknot) using two negative crossing changes, but an analysis
of signatures and a minor generalization of the results of \cite{cl} shows that
it cannot be converted into a slice knot (or a knot with signature 0) using two
positive generalized crossing changes. 

Related to this discussion we have the follow result.  Its proof is a bit
technical to include here and will be described in detail elsewhere.

\begin{theorem} A knot $K$ with 3--sphere genus $g(K)$ can be converted into
the unknot using
$2g(K)$ generalized crossing changes. \end{theorem}
 
\bf Addendum\qua\rm (December 15, 2002)\qua It has been pointed out to the
author that results of Murakami and Yasuhara ({\it Four-genus and
four-dimensional clasp number of a knot}, Proc. Amer. Math. Soc. 128
(2000), no. 12, 3693--3699) imply that $g_s(8_{16}) = 1$ but
$u_s(8_{16}) =2$. The methods used there are different from those of
this paper.


\Addressesr

\end{document}